\documentclass[12pt, onesided, reqno]{amsart}

\usepackage{amssymb}
\usepackage{amsmath}
\usepackage{color}
\usepackage{enumerate}

\evensidemargin\oddsidemargin

\begin{document}

\setcounter{page}{1}

\newtheorem{PROP}{Proposition}
\newtheorem{REM}{Remark}
\newtheorem{LEM}{Lemma}
\newtheorem{THEA}{Theorem A\!\!}
\renewcommand{\theTHEA}{}
\newtheorem{THEB}{Theorem B\!\!}
\renewcommand{\theTHEB}{}
\newtheorem{CORA}{Corollary 1\!\!}
\renewcommand{\theCORA}{}
\newtheorem{CORB}{Corollary 2\!\!}
\renewcommand{\theCORB}{}

\newtheorem{theorem}{Theorem}
\newtheorem{proposition}[theorem]{Proposition}
\newtheorem{corollary}[theorem]{Corollary}
\newtheorem{lemma}[theorem]{Lemma}
\newtheorem{assumption}[theorem]{Assumption}

\newtheorem{definition}[theorem]{Definition}
\newtheorem{hypothesis}[theorem]{Hypothesis}

\theoremstyle{definition}
\newtheorem{example}[theorem]{Example}
\newtheorem{remark}[theorem]{Remark}
\newtheorem{question}[theorem]{Question}

\newcommand{\eqnsection}{
\renewcommand{\theequation}{\thesection.\arabic{equation}}
    \makeatletter
    \csname  @addtoreset\endcsname{equation}{section}
    \makeatother}
\eqnsection

\def\a{\alpha}
\def\b{\beta}
\def\B{{\bf B}} 
\def\car{c_{\a, \rho}}
\def\cah{c_{\a, \hr}}
\def\cua{c_{\a, u}}
\def\cL{{\mathcal{L}}} 
\def\Ea{E_\a}
\def\eps{\varepsilon}
\def\g{{\gamma}} 
\def\hr{{\hat{\rho}}}
\def\hL{{\hat{L}}}
\def\Ga{{\Gamma}} 
\def\i{{\rm i}}
\def\K{{\bf K}}
\def\Ka{{\bf K}_\a}
\def\kar{{\kappa_{\a, \rho}}}
\def\L{{\bf L}}
\def\lbd{\lambda}
\def\lcr{\left[}
\def\lpa{\left(}
\def\lva{\left|}
\def\otheta{\overline{\theta}}
\def\rpa{\right)}
\def\rcr{\right]}
\def\rva{\right|}
\def\T{{\bf T}}
\def\M{{\mathcal M}}
\def\X{{\bf X}}
\def\Ton{{T_0^{(n)}}}
\def\Lton{{\vert L_{T_0^{(n)}}\vert}}
\def\Un{{\bf 1}}
\def\ZZ{{\bf Z}}
\def\CC{{\bf C}}
\def\GG{{\bf \Ga}}
\def\BB{{\bf B}}
\def\car{c_{\a,\rho}}
\def\sar{s_{\a,\rho}}
\def\pbxy{\Pb_{(x,y)}}
\def\foxy{f^0_{x,y}}

\def\E{\mathbb{E}}
\def\Pb{\mathbb{P}}
\def\esp{\mathbb{E}}
\def\Z{\mathbb{Z}}
\def\N{\mathbb{N}}
\def\Q{\mathbb{Q}}
\def\rl{\mathbb{R}}
\def\P{\mathcal{P}}
\def\pb{\mathbb{P}}
\def\C{\mathbb{C}}
\def\cC{\mathcal{C}}
\def\F{\mathcal{F}}
\def\S{\mathcal{S}}
\def\W{\mathcal{W}}
\def\L{\mathcal{L}}
\def\G{\mathcal{G}}

\newcommand*\pFqskip{8mu}
\catcode`,\active
\newcommand*\pFq{\begingroup
        \catcode`\,\active
        \def ,{\mskip\pFqskip\relax}%
        \dopFq
}
\catcode`\,12
\def\dopFq#1#2#3#4#5{%
        {}_{#1}F_{#2}\biggl[\genfrac..{0pt}{}{#3}{#4};#5\biggr]%
        \endgroup
}

\newcommand{\equi}{\mathop{\sim}\limits}
\def\={{\;\mathop{=}\limits^{\text{(law)}}\;}}
\def\d{{\;\mathop{=}\limits^{\text{(1.d)}}\;}}
\def\st{{\;\mathop{\geq}\limits^{\text{(st)}}\;}}

\def\cas{\stackrel{a.s.}{\longrightarrow}}
\def\claw{\stackrel{d}{\longrightarrow}}
\def\elaw{\stackrel{d}{=}}
\def\qed{\hfill$\square$}
                  
\title[Harmonic measure of stable processes]
       {On the harmonic measure of stable processes}
       
\author[Christophe Profeta]{Christophe Profeta}

\address{Laboratoire de math\'ematiques et mod\'elisation d'Evry, Universit\'e d'Evry-Val d'Essonne, F-91037 Evry Cedex. {\em Email} : {\tt christophe.profeta@univ-evry.fr}}

\author[Thomas Simon]{Thomas Simon}

\address{Laboratoire Paul Painlev\'e, Universit\'e Lille 1, F-59655 Villeneuve d'Ascq Cedex and Laboratoire de physique th\'eorique et mod\`eles statistiques, Universit\'e  Paris-Sud, F-91405 Orsay Cedex. {\em Email} : {\tt simon@math.univ-lille1.fr}}



\begin{abstract} Using three hypergeometric identities, we evaluate the harmonic measure of a finite interval and of its complementary for a strictly stable real L\'evy process. This gives a simple and unified proof of several results in the literature, old and recent. We also provide a full description of the corresponding Green functions. As a by-product, we compute the hitting probabilities of points and describe the non-negative harmonic functions for the stable process killed outside a finite interval.
\end{abstract}

\maketitle
 
\section{Introduction and statement of the results}

Let $L = \{L_t, \, t\ge 0\}$ be a real strictly $\a-$stable L\'evy process, with characteristic exponent
\begin{equation}
\label{Norm}
\Psi(\lbd)\; =\;\log(\esp[e^{\i \lbd L_1}])\; =\;-\; (\i \lbd)^\a e^{-\i\pi\a\rho\, {\rm sgn}(\lbd)}, \qquad \lbd\in\rl.
\end{equation}
Above, $\a\in (0,2]$ is the self-similarity parameter and $\rho = \pb[L_1\ge 0]$ is the positivity parameter. Recall that when $\a = 2,$ one has $\rho = 1/2$ and $\Psi(\lbd) = -\lbd^2,$ so that $L$ is a Brownian motion. When $\a = 1,$ one has $\rho\in (0,1)$ and $L$ is a Cauchy process with a linear drift. When $\a\in (0,1)\cup(1,2)$ the characteristic exponent reads
$$\Psi(\lbd) \; =\; \; -\; \kappa_{\a,\rho}\vert\lbd\vert^\a (1 - \i\b\tan(\pi\a/2)\,{\rm sgn}(\lbd)),$$
where $\b\in[-1,1]$ is the asymmetry parameter, whose connection with the positivity parameter is given by Zolotarev's formula:
$$\rho \; =\; \frac{1}{2} \,+ \,\frac{1}{\pi\a} \arctan(\b\tan(\pi\a/2)),$$
and $\kappa_{\a,\rho} = \cos(\pi\a(\rho -1/2)) > 0$ is a scaling constant. We refer e.g. to Chapter VIII in \cite{B} for more details on this parametrization. One has $\rho \in [0,1]$ if $\a < 1$ and $\rho\in[1-1/\a, 1/\a]$ if $\a > 1.$ When $\a > 1, \rho = 1/\a$ or $\a < 1, \rho =0,$ the process $L$ has no positive jumps, whereas it has no negative jumps when $\a > 1, \rho = 1-1/\a$ or $\a < 1, \rho =1.$

Set $\hL = -L$ for the dual process and $\hr = 1-\rho$ for its positivity parameter. Throughout, it will be implicitly assumed that all quantities enhanced with a hat refer to the same quantities for the dual process, that is with $\rho$ and $\hr$ switched. We denote by $\pb_x$ the law of $L$ starting from $x\in\rl.$ Introduce the harmonic measures
$$H_x(dy)\; =\; \pb_x [L_T \in dy, \; T< \infty]\qquad \text{and}\qquad H_x^*(dy)\; =\; \pb_x [L_{T^*} \in dy, \; T^*< \infty],$$
where $T= \inf\{ t > 0, \; \vert L_t \vert > 1\}$ and $T^* = \inf\{t>0,\; \vert L_t \vert < 1\}.$ Observe that by spatial homogeneity and the scaling relationship
\begin{equation}
\label{Sca}
\lpa \{kL_t, \, t \ge 0\}, \pb_x\rpa\; \elaw\; \lpa \{L_{k^\a t}, \, t \ge 0\}, \pb_{kx}\rpa, \qquad k >0,
\end{equation}
we can deduce from $H_x$ the expression of the harmonic measure of the complementary of any closed bounded interval, whereas the knowledge of $H_x^*$ gives that of the harmonic measure of any open bounded interval. Introduce the following notation
$$x_+\; =\; \max(x,0), \qquad c_{\a,\rho}\; =\; \frac{\sin (\pi\a\rho)}{\pi} \qquad\mbox{and}\qquad \psi_{\a,\rho}(t) \;= \; (t-1)^{\a\hr-1}(t+1)^{\a\rho-1}.$$
In the remainder of this section it will be implicitly assumed that $L$ has jumps of both signs. The corresponding results where $L$ has one-sided jumps, which are simpler, will be briefly described in the last section.

\begin{THEA} (a) For any $x \in(-1, 1),$ the measure $H_x(dy)$ has density 
$$h(x,y) \; =\; c_{\a,\rho} \,(1+x)^{\a\hr}(1-x)^{\a\rho}(1+y)^{-\a\hr}(y-1)^{-\a\rho}(y-x)^{-1}$$

\smallskip

\noindent
if $y > 1$ and $h(x,y) = {\hat h(-x,-y)}$ if $y <-1.$\\

(b) For any $x\in [-1,1]^c,$ the measure $H_x^*(dy)$ has density 
$$h^*(x,y) \; =\; c_{\a,\hr}\, (1+y)^{-\a\rho}(1-y)^{-\a\hr}\lpa (x+1)^{\a\rho}(x-1)^{\a\hr}(x-y)^{-1}  - (\a-1)_+ \!\int_1^x \!\psi_{\a,\rho}(t)\, dt\rpa$$
if $x > 1,$ and $h^*(x,y) = {\hat h^*(-x,-y)}$ if $x <-1.$

\end{THEA}

In the symmetric case, these computations date back to \cite{BGR} - see Theorems A, B and C therein. Notice that the results of \cite{BGR}, which rely on  Kelvin's transformation and the principle of unicity of potentials, deal with the more general rotation invariant stable processes on Euclidean space. In the general case, Part (a) of the above theorem was proved in Theorem 1 of \cite{Ro}, whereas Part (b) was recently obtained in Theorem 1 of \cite{KPW}. Both methods used in \cite{Ro} (coupled integral equations) and in \cite{KPW} (Lamperti's representation and the Wiener-Hopf factorization) are complicated. In this paper we show that the original method of \cite{BGR} works in the asymmetric case as well, thanks to elementary considerations on the hypergeometric function
$$\pFq{2}{1}{a,b}{c}{z}\; =\; \sum_{n\ge 0} \frac{(a)_n(b)_n}{(c)_n}\, z^n.$$ 
More precisely, we use three basic identities for the latter function, due respectively to Euler, Pfaff and Gauss, allowing to perform a simple potential analysis of the function
\begin{equation}
\label{Pho}
\varphi(t)\; =\; (1-t)^{-\a\rho}(1+t)^{-\a\hr}
\end{equation}
and to obtain the required generalization of the key Lemma 3.1 in \cite{BGR}. \\

Define next the killed potential measures
$$G_x(dy)\; =\; \esp_x\lcr \int_0^T \Un_{\{L_t\in dy\}} \,dt\rcr \qquad \text{ and }\qquad G_x^*(dy)\; =\; \esp_x\lcr \int_0^{T^*} \Un_{\{L_t\in dy\}} \,dt\rcr.$$
It is easy to see from the absolute continuity of the two killed semi-groups with respect to the original stable semigroup, that both these measures are absolutely continuous. We denote by $g(x,y)$ and $g^*(x,y)$ their respective densities on $(-1,1)$ and $[-1,1]^c,$ the so-called Green functions. These functions are of central interest because they allow to invert the stable infinitesimal generator on $(-1,1)$ and on $[-1,1]^c$ - see e.g. Formula (1.42) in \cite{BBKRSV} in the symmetric case. Observe that they are related to the harmonic measure and to the density of the L\'evy measure of $L$:
\begin{equation}
\label{Nu}
\nu(y)\; =\; \Ga(\a +1) \vert y\vert^{-\a-1}\lpa c_{\a,\rho}\Un_{\{ y > 0\}} \, +\, c_{\a,\hr}\Un_{\{ y < 0\}}\rpa,
\end{equation}
through the integral formul\ae
$$h(x,y)\; =\; \int_{(-1,1)} \! g(x,v) \nu(y-v)\, dv\qquad \mbox{resp.}\qquad h^*(x,y)\; =\; \int_{[-1,1]^c} \!g^*(x,v) \nu(y-v) \,dv$$
for all $x\in (-1,1)$ and $y \in [-1,1]^c$ resp. for all $x\in [-1,1]^c$ and $y \in (-1,1),$ which are both instances of a general formula by Ikeda-Watanabe - see Theorem 1 in \cite{IW}. For this reason, the density of the harmonic measure coincides with that of the Poisson kernel - see \cite{BBKRSV} pp. 16-17. The closed expression of the Poisson kernel and the Green function for $(-1,1)$ in the symmetric case, and more generally for the open unit ball in the rotation invariant case, are classic results dating back to Riesz \cite{Ri, Ri1}. We refer to \cite{BBKRSV} pp. 18-19 for more details and references, and to the whole monograph \cite{BBKRSV} for several extensions, all in the rotation invariant framework. 
    
\begin{THEB} Set $z = z(x,y) = \lva \frac{1-xy}{x-y} \rva$ for every $x\neq y.$ 

\begin{enumerate}[$(a)$]
\item For every $x \in (-1 ,1),$ one has 
$$g(x,y) \; =\; \frac{1}{\Ga(\a\rho)\Ga(\a\hr)}\, \lpa \frac{y-x}{2}\rpa^{\a-1}\int_1^z \psi_{\a,\rho}(t)\, dt$$
if $y \in(x,1),$ and $g(x,y) = {\hat g}(y,x)$ if $y\in (-1,x).$ \\

\item For every $x > 1,$ one has 
$$g^*(x,y)\; =\;  \frac{2^{1-\a}}{\Ga(\a\rho)\Ga(\a\hr)}\, \lpa (y-x)^{\a-1}\int_1^z\psi_{\a,\rho}(t) dt \, -\, (\a-1)_+\int_1^x \psi_{\a,\rho}(t) dt\int_1^{y } \psi_{\a,\hr}(t)\, dt\rpa$$
if $y \in(x,\infty), g(x,y) = {\hat g}(y,x)$ if $y\in (1,x),$ and
$$g^*(x,y)\; =\;  \frac{c_{\a,\hr}\, 2^{1-\a}}{c_{\a,\rho}\,\Ga(\a\rho)\Ga(\a\hr)}\lpa (x-y)^{\a-1}\!\int_1^z\!\psi_{\a,\rho}(t) dt  - (\a-1)_+\!\int_1^x\! \psi_{\a,\rho}(t) \,dt\int_1^{\vert y\vert } \!\!\psi_{\a,\rho}(t)\, dt\rpa$$
if $y <-1.$

\end{enumerate}

\end{THEB}

Observe that in Part (b) of the above result, the condition $x> 1$ is no restriction since by duality we have $g^*(x,y) = {\hat g^*}(-x,-y)$ for every $x <-1$ and $y\in [-1,1]^c.$ Part (a) was obtained as Corollary 4 of \cite{BGR} in the symmetric case, and as Theorem 1 of \cite{KW} in the general case. Part (b) was proved as Theorem 4 in \cite{KPW}, in the only cases $\a \le 1$ and $x,y > 1.$ The methods of \cite{KPW, KW}, relying on the Lamperti transformation and an analysis of the reflected process, are complicated. In this paper, we observe that all formul\ae \, of Theorem B can be quickly obtained from the D\'esir\'e Andr\'e equation and one of the two simple lemmas leading to the proof of Theorem A. \\

The explicit knowledge of the Green function has a number of classical consequences. In this paper we will focus on two of them. The first one deals, in the relevant case $ \a > 1,$ with the hitting probability $\rho(x,y)  = \pb_x [T_y < T],$ where $T_{y} = \inf\{ t > 0, \, L_t =y\}.$
 
\begin{CORA} Assume $\a > 1$ and set $z = \lva \frac{1-xy}{x-y} \rva.$ For every $x, y \in (-1,1),$ one has
$$\rho(x,y) \; =\; (\a -1)\, \lpa \frac{x-y}{1-y^2}\rpa^{\a-1} \int_1^z \psi_{\a, \hr} (t)\, dt$$
if $x > y,$ and $\rho(x,y)= {\hat \rho}(-x,-y)$ if $x < y.$ 
\end{CORA}
Observe that the above formula extends by continuity on the diagonal, with the expected property that $\pb_x[T_x < T] = 1.$ Of course, this follows from the fact that $\{x\}$ is regular for $x$ in the case $\a > 1.$ When $\a \to 2,$ Corollary 1 amounts to the very standard Brownian formula
$$\pb_x[T_y < T_1] \; =\; \frac{1-x}{1-y}\cdot$$
By the Markov property, one can deduce from Corollary 1 the harmonic measure of the set $\{y\}\, \cup \,[-1,1]^c.$ Using one of our three hypergeometric identities, it is also possible to derive the asymptotic behaviour of $\pb_x [T_y<T]$ when $x\to y,$ which is fractional. Last, by spatial homogeneity and scaling, we can quickly recover the statement of Theorem 1.5 in \cite{KPW}. See Remark 6 below for more detail.  \\

We next consider non-negative harmonic functions on $(-1,1),$ which are the non-negative solutions to 
$$\cL_{\a, \rho} u \; \equiv\; 0$$
on $(-1,1),$ where $\cL_{\a,\rho}$ is the infinitesimal generator of $L$. As in the Brownian case, an equivalent characterization - see e.g. \cite{BBKRSV} p. 20 in the symmetric case - is the mean-value property, which reads $\esp_x[u(L_{\tau_U})] = u(x)$ for every open set $U$ whose closure belongs to $(-1,1),$ where $\tau_U = \inf\{ t > 0,\, L_t\not\in U\}.$

\begin{CORB} The non-negative harmonic functions on $(-1,1)$ which vanish on $[-1,1]^c$ are of the type
$$x\; \mapsto\;\lbd (1-x)^{\a\rho}(1+x)^{\a\hr -1} \; +\; \mu (1+x)^{\a\hr} (1-x)^{\a\rho-1}$$
with $\lbd, \mu \ge 0.$
\end{CORB}

This result might be already known - compare e.g. with Theorem 10 p.569 in \cite{Si}, although we could not find it written down explicitly in the literature. Recall that in order to obtain all non-negative harmonic functions on $(-1,1)$, one needs  - see e.g. Theorem 2.6 in \cite{BBKRSV} in the symmetric case - to add to the above functions the integral of the Poisson kernel $h(x,y)$ along some suitably integrable measure  on $[-1,1]^c.$ \\

Both Corollaries 1 and 2 could be obtained for the process killed inside the interval $(-1,1),$ with analogous computations relying on Part (b) of Theorem B. But the formul\ae \, have a rather lengthy aspect, so that we prefer leaving them to the interested reader. The remainder of the paper is as follows. In the three next sections we prove Theorem A, Theorem B, and the two Corollaries. In the last section we gather, for the sake of completeness, the corresponding formul\ae\, in the cases of semi-finite intervals and of one-sided jumps.
  
\section{Proof of Theorem A} 

As mentioned in the introduction, the argument hinges upon three classical hypergeometric identities, to be found in Theorem 2.2.1, Formula (2.2.6) and Formula (2.3.12) of \cite{AAR}, and which will be henceforth referred to as Euler, Pfaff and Gauss\footnote{among of course many others. This one is a simple consequence of the two-dimensional structure of the space of solutions to the hypergeometric equation. Notice that it can also be obtained by Mellin-Barnes inversion. See the end of the article {\em Calculs asymptotiques} in Encyclopedia Universalis.} formula respectively. 

\subsection{Proof of Part (b)}

\subsubsection{The case $\a < 1$.} We reason along the same lines as in Theorem B of \cite{BGR}. Set $p_t(x)$ for the transition density of $L.$ The following computation, left to the reader, is a well-known consequence of (\ref{Norm}), Fourier inversion and the Fresnel integral: one has
$$\int_0^\infty p_t(z)\, dt\; =\; \Ga(1-\a)\, c_{\a,\rho} \, z^{\a -1}$$
for every $z > 0.$ Observe that by duality, one also has
$$\int_0^\infty p_t(z)\, dt\; =\; \int_0^\infty {\hat p}_t(-z)\, dt\; =\; \Ga(1-\a) \, c_{\a,\hr} \,  \vert z\vert^{\a -1}$$
for every $z < 0.$ Applying the D\'esir\'e Andr\'e equation (2.1) in \cite{BGR} and letting $s\to 0$ therein shows that
\begin{equation}
\label{Abel1}
\int_{-1}^1 u(t,y) \, H_x^*(dt) \; =\; \cah\, \vert x-y\vert^{\a -1} \end{equation}
for every $x > 1$ and $y\in (-1,1),$ where we have set 
$$u(t,y)\; =\; \lpa \car\, \Un_{\{y > t\}} + \cah \,\Un_{\{y < t\}} \rpa \vert t-y\vert^{\a -1}.$$ 
In the symmetric case, this Abelian integral equation with constant boundary terms is solved in Section 3 of \cite{BGR}, following the method of \cite{Ri}. See also \cite{C} for the original solution, with a more general term on the left-hand side. After proving the following lemma, which remains valid for $\a\in (1,2),$ we will see that the pole-seeking method of \cite{Ri} applies in the asymmetric case as well.
 
\begin{LEM} The unique positive measure on $(-1,1)$ satisfying
\begin{equation}
\label{Abel0}
\int_{-1}^1 {\hat u}(t,y) \, \mu(dt)\; = \; 1, \qquad y\in (-1,1) 
\end{equation}
has the density $\varphi(t)$ given in {\em (\ref{Pho}).}
\end{LEM}

\proof The fact that there is a unique measure solution of (\ref{Abel0}) is a standard fact in potential theory - see e.g. Theorem 1 in \cite{P1} or Proposition VI.1.15 in \cite{BG}. In our concrete context, this unicity can also be obtained by a straightforward adaptation of Lemma 4.1 in \cite{BGR}. To show the lemma, we compute by a change of variable
\begin{eqnarray*}
\int_{-1}^1 {\hat u}(t,y) (1-t)^{-\a\rho}(1+t)^{-\a\hr}\, dt & = & \cah\int_{0}^{1+y} t^{\a-1}(1-y+t)^{-\a\rho}(1+y-t)^{-\a\hr}\, dt\\
& & +\;\;  \car \int_0^{1-y} t^{\a-1}(1-y-t)^{-\a\rho}(1+y+t)^{-\a\hr}\,dt.
\end{eqnarray*}
Using two further changes of variable, the Euler formula, and the complement formula for the Gamma function, we transform the expression on the right-hand side into
$$\frac{\Ga(\a)}{\Ga(\a\hr)\Ga(1+\a\rho)} \lpa\frac{1+y}{1-y}\rpa^{\a\rho} \lpa \pFq{2}{1}{\a\rho,\a}{1+\a\rho}{\frac{y+1}{y-1}}\; +\; \frac{\rho}{\hr} \lpa\frac{1-y}{1+y}\rpa^\a \pFq{2}{1}{\a\hr,\a}{1+\a\hr}{\frac{y-1}{y+1}}\rpa,$$
and then, using the notation
$$z\; =\; \frac{y+1}{y-1},$$ 
into
$$\frac{\Ga(\a) (-z)^{\a\rho}}{\Ga(\a\hr)\Ga(1+\a\rho)} \lpa \pFq{2}{1}{\a, \a\rho}{1+\a\rho}{z}\; +\; \frac{\rho (-z)^{-\a}}{\hr} \,\pFq{2}{1}{\a,\a\hr}{1+\a\hr}{\frac{1}{z}}\rpa\; =\; 1,$$
\smallskip

\noindent
where the last equality follows from the Gauss formula.

\endproof

\begin{remark}The solution to (\ref{Abel0}) in the symmetric case was obtained in Lemma 3.1 of \cite{BGR}, via a reflection argument. Alternatively, the non-symmetric solution can be deduced in a constructive way, following the approach of \cite{Ri} pp. 41-42 or that of \cite{C}. Observe that the above argument is significantly shorter than in these three references. 
\end{remark}

We can now finish the proof. Introduce the changes of variables
\begin{equation}\label{cov}
t \; =\; x \; + \; \frac{1-x^2}{x-s}\qquad\mbox{and}\qquad y \; =\; x \; + \; \frac{1-x^2}{x-z},
\end{equation}
and observe that they map $(-1,1)$ onto $(-1,1),$ in a decreasing way. Plugging these changes of variables into (\ref{Abel0}) implies after some computation that
$$(x+1)^{\a\rho}(x-1)^{\a\hr} \, \int_{-1}^1  (1+s)^{-\a\rho}(1-s)^{-\a\hr}(x-s)^{-1}\, u(s,z)\, ds\; =\; \vert x-z\vert^{\a -1}$$ 
for every $x > 1$ and $z\in (-1,1).$ Multiplying both sides by $\cah$ shows the required solution to (\ref{Abel1}), which is unique by Lemma 1 and the changes of variables (\ref{cov}). 

\qed

\begin{remark}
In the following, we shall make a repeated use of the changes of variable (\ref{cov}), which may be written formally :
$$|1+x|^{\a\rho}|1-x|^{\a\hr}\, \int  |y-t| ^{\a-1} \frac{|1+t|^{-\a\rho}|t-1|^{-\a\hr}}{|x-t|} \,dt \;=\; 
|x-y|^{\a -1} \int  |z-s|^{\a -1} \varphi(s) ds.$$
The interest of this change of variable is to transform an Abelian integral with two inside parameters into an integral of the hypergeometric type, with one parameter inside.
\end{remark}

\medskip

\subsubsection{The case $\a > 1$.} We follow the method of Theorem C in \cite{BGR}. Recall that since $L$ a.s. hits points in finite time, the measure $H^*_x(dt)$ has total mass one. We will need the evaluation
$$\int_0^\infty (p_t(z) - p_t(0))\, dt\; =\; \Ga(1-\a) \,\car\, z^{\a -1}$$
for every $z >0,$ which is easy and classical - see the introduction of \cite{P1}. This implies
$$\lpa\int_0^\infty e^{-st} p_t(z) \, dt \; -\; p_1(0)\Ga(1- 1/\a) s^{\frac{1}{\a} -1}\rpa\; \downarrow\; \Ga(1-\a) \,\car\, z^{\a -1}$$ 
as $s\to 0,$ for every $z >0.$ Proceeding as in \cite{BGR} pp. 544-545 shows that
\begin{equation}
\label{Abel2}
\cah\, \vert x-y\vert^{\a -1} \; =\; \int_{-1}^1  u(t,y)\, H^*_x(dt)  \; +\; \kappa_{\a,\rho}^* (x)
\end{equation}
for every $x > 1$ and $y\in (-1,1),$ where 
$$\kappa_{\a,\rho}^*(x)\;=\; \frac{p_1(0) \Gamma(1-1/\alpha)}{\Gamma(1-\alpha)}\, \times\, \lim_{\lambda \rightarrow 0} \lambda^{1/\alpha-1} \lpa \E_x\left[e^{-\lambda T^*}\right] - 1 \rpa$$
is a non-negative function which will be determined in the same way as in (4.1) of \cite{BGR}. Multiplying both sides of (\ref{Abel2}) by $\varphi(y)$ and integrating on $(-1,1)$ shows by Lemma 1 that
$$\kappa_{\a,\rho}^* (x) \;= \; \lpa \int_{-1}^1 \varphi(y)\, dy\rpa^{-1}\lpa\cah \int_{-1}^1 (x-y)^{\a -1}\varphi(y)\, dy\; -\; 1\rpa$$
for every $x > 1.$ The next lemma, generalizing the second part of Lemma 3.1 in \cite{BGR}, allows to compute the right-hand side.

\begin{LEM}\label{lem:phi} With the above notation, one has
$$\cah\, \int_{-1}^1 (x-y)^{\a -1}\varphi(y)\, dy\; =\; 1\; -\;\frac{\Ga(1-\a\rho)2^{1-\alpha}}{\Ga(\a\hr)\Ga(1-\a)}\;  \int_1^x \psi_{\a,\rho}(t)\, dt$$
for every $x > 1.$
\end{LEM}

\proof As in Lemma 1, a change of variable and the Euler formula show first that
\begin{eqnarray*}
\sin (\pi\a\hr)\int_{-1}^1 (x-y)^{\a -1}\varphi(y)\, dy & = & \frac{\Ga(1-\a\rho)}{\Ga(\a\hr)\Ga(2-\a)}\lpa \frac{x+1}{2}\rpa^{\a-1} \!\!\pFq{2}{1}{1-\a,1-\a\hr}{2-\a}{\frac{2}{x+1}}\\
& = & \frac{\Ga(1-\a\rho)}{\Ga(\a\hr)\Ga(2-\a)}\lpa \frac{x-1}{2}\rpa^{\a-1} \!\!\pFq{2}{1}{1-\a,1-\a\rho}{2-\a}{\frac{2}{1-x}},
\end{eqnarray*}
where the second equality follows from the Pfaff formula. Using now the Gauss formula, we next transform 
$$(-z)^{\a-1} \pFq{2}{1}{1-\a,1-\a\rho}{2-\a}{\frac{1}{z}}\; =\;  
\frac{\Ga(\a\hr)\Ga(2-\a)}{\Ga(1-\a\rho)}\; +\; \frac{(\a-1)}{\a\hr}(-z)^{\a\hr} \pFq{2}{1}{1-\a\rho,\a\hr}{1+\a\hr}{z}$$
with the notation $z= (1-x)/2.$ Putting everything together and applying again the Euler  formula completes the proof.

\endproof

We can now finish the proof of the case $\a > 1$. From Lemma 2 and an easy computation, we first deduce 
$$\kappa_{\a,\rho}^*(x)\; =\; \cah\, 2^{\a -1}\,\frac{1}{\a\hr}\lpa\frac{x-1}{2}\rpa^{\a\hr }\!\! \pFq{2}{1}{1-\a\rho,\a\hr}{1+\a\hr}{\frac{1-x}{2}}.$$
Coming back to (\ref{Abel2}) and reasoning as in \cite{BGR} p.552, we finally see from Lemma 1 that $H_x^* (dy)$ has density
$$\cah\, (x+1)^{\a\rho}(x-1)^{\a\hr}(1+y)^{-\a\rho}(1-y)^{-\a\hr}(x-y)^{-1} \; -\; \kappa_{\a,\rho}^*(x){\hat \varphi}(y).$$
To conclude the proof, it suffices to observe by the Euler formula and a change of variable that 
\begin{eqnarray*}
\kappa_{\a,\rho}^*(x){\hat \varphi}(y) & =& \cah\, (1+y)^{-\a\rho}(1-y)^{-\a\hr}(\a-1)\, 2^{\a -1}\,\frac{1}{\a\hr}\lpa\frac{x-1}{2}\rpa^{\a\hr }\!\! \pFq{2}{1}{1-\a\rho,\a\hr}{1+\a\hr}{\frac{1-x}{2}}\\
& = & \cah\, (1+y)^{-\a\rho}(1-y)^{-\a\hr}\, (\a-1) \int_1^x \psi_{\a,\rho}(t)\, dt.
\end{eqnarray*}
\qed

\begin{remark}
Since $\kappa_{\alpha, \rho}^* (x)$ is finite and positive, we can deduce from Karamata's Tauberian theorem that  
$$\Pb_x[T^*>t]\;\sim\;  - \frac{\Gamma(1-\a)\sin(\pi/\a)}{\pi p_1(0)} \, \kappa_{\a, \rho}^*(x) \,t^{1/\a-1}\qquad \mbox{as} \;t\to +\infty.$$
This asymptotic is given in Corollary 3 of \cite{BGR} in the symmetric case, and in Theorem 2 of \cite{P1} in the asymmetric case, with a more general formulation. Notice that $T^*$ has infinite expectation.

\end{remark}

\medskip

\subsubsection{The case $\a =1$.} This case is known to be more subtle from the computational point of view, because it involves logarithmic kernels. The transition density of $L_t$ is
$$p_t (x) \; =\; \frac{c_{1,\rho}\,t}{t^2 + 2tx \cos\pi\rho + x^2}\cdot$$
The process $L$ does not hit points a.s. but it is recurrent, so that $H_x^* (dt)$ has total mass one. After some computation, one finds
$$\int_0^\infty (p_t(1) - p_t (x))\, dt\; =\; c_{1,\rho}\,\log\vert x\vert .$$
See also \cite{P1} p.391. With this formula, it is possible to finish the proof as in the case $\a > 1,$ but the computations are lenghty and we hence prefer to invoke a simple argument relying on the Skorokhod topology. Fix $\rho \in (-1,1)$ and let $\a\downarrow 1.$ It follows from Corollary VII.3.6 in \cite{JS} that
$$\cL (L^{\a, \rho})\; \Rightarrow \; \cL (L^{1,\rho})$$
with obvious notation for $L^{1, \rho}$ and $L^{\a, \rho}$, and where $\Rightarrow$ means weak convergence in the classical Skorokhod space. Using Remark VI.3.8 and Proposition VI.2.12 in \cite{JS}, it is then easy to deduce that
$$L_T^{\a,\rho} \; \claw\; L_T^{1,\rho}.$$
The conclusion follows from pointwise convergence of the densities $h^*(x,y)$ as $\a\downarrow 1,$ and Scheff\'e's lemma.

\qed

\begin{remark} The above argument relying on a.s. continuity for the Skorokhod topology will be used repeatedly in the sequel, under the denomination Skorokhod continuity argument.
\end{remark}

\smallskip

\subsection{Proof of Part (a)} By the Skorokhod continuity argument, it is enough to consider the case $\a\neq 1.$ Fixing $x \in (-1,1)$ and proceeding as in \cite{BGR} pp. 544-545, the harmonic measure $H_x(dt)$ is seen to be the unique solution of the equation
\begin{equation}
\label{Abel3}
u(x,y) \; =\; \int_{(-1,1)^c} u(t,y) \, H_x(dt)
\end{equation}
for every $y \in [-1,1]^c.$ In the case $\a < 1,$ this is indeed an immediate consequence of the Markov property, leading to the corresponding equation (\ref{Abel1}). And in the case $\alpha>1$, the well-known fact - see Lemma 4.1 in \cite{Ta} -  that the tail distribution of $T$ is exponentially small at infinity implies that the perturbative term $\kappa_{\alpha,\rho}$ is zero in the corresponding equation (\ref{Abel2}). Define 
$$\mu_x(dt) \; = \;
\begin{cases}
\car\, \hat{\varphi}(t)\, dt &\qquad \text{if }t\leq x,\\
\cah\, \hat{\varphi}(t)\, dt &\qquad \text{if }t> x.\\
\end{cases}
$$
We shall deal with the two cases $y>1$ and $y<-1$ separately. \\

\noindent
$(i)$ Let $\nu\in (-1, x)$.  Applying Lemma 1 with $\rho$ and $\hr$ interchanged, we get
\begin{align*}
&\car\int_{-1}^\nu |\nu-t|^{\alpha-1}\mu_x(dt)\, + \,\cah \int_{\nu}^x |\nu-t|^{\alpha-1}\mu_x(dt)\,+\,\car\int_{x}^1 |\nu-t|^{\alpha-1}\mu_x(dt)\\
&\qquad =\;\car \left(\int_{-1}^\nu  \car |\nu-t|^{\alpha-1} \hat{\varphi}(t) dt \, + \,\int_\nu^{1} \cah |\nu-t|^{\alpha-1} \hat{\varphi}(t) dt\right)\;=\;\car .
\end{align*}
The changes of variable (\ref{cov}) implies after some rearrangement
$$\int_{[-1,1]^c} \vert y-t\vert^{\alpha-1}\, H_x(dt) \; =\;\car(y-x)^{\a-1}$$
for every $y > 1$ with the required expression for $H_x(dt),$ which is Equation (\ref{Abel3}). \\

\noindent
$(ii)$ Take now $\nu \in (x, 1)$. Applying again Lemma 1, we have
\begin{align*}
&\cah\int_{-1}^x |\nu-t|^{\alpha-1}\mu_x(dt)\, + \,\car \int^{\nu}_x |\nu-t|^{\alpha-1}\mu_x(dt)\,+\,\car\int_{\nu}^1 |\nu-t|^{\alpha-1}\mu_x(dt)\\
&\qquad =\;\cah \left(\int_{-1}^\nu  \car |\nu-t|^{\alpha-1} \hat{\varphi}(t) dt \, + \,\int_\nu^{1} \cah |\nu-t|^{\alpha-1} \hat{\varphi}(t) dt\right)\;=\;\cah.
\end{align*}
The same changes of variables (\ref{cov}) gives
$$\int_{[-1,1]^c} \vert y-t\vert^{\alpha-1}\, H_x(dt) \; =\;\cah(x-y)^{\a-1}$$
for every $y < -1,$ which is again Equation (\ref{Abel3}).

\qed

\begin{remark} (a) The behaviour at infinity of the distribution function of $T$ is more mysterious than that of $T^*$. In the non-subordinator case it is known - see Proposition VIII.3 in \cite{B} - that there exists $\kappa_x$ positive and finite such that
$$-\log \pb_x [T > t]\; \sim\; -\kappa_x t\qquad \mbox{as $t\to +\infty,$}$$
but the exact value of $\kappa_x$ is unknown except in the completely asymmetric case - see \cite{B1}. We refer to Chapter 4 in \cite{BBKRSV} for more on this topic in the rotation invariant case. Notice that the result of Theorem B (a) allows to compute the expectation of $T$:
$$\E_x[T]\; =\; \int_{-1}^1 g(x,y)\, dy\; =\; \frac{(1-x)^{\a \rho}(1+x)^{\a\hr}}{\Gamma(\a+1)}\cdot$$

\medskip

(b) With our computations, we can also check the values of the total masses $H_x(-1,1)^c$ and $H^*_x(-1,1).$ On the one hand, Lemma 1 and the change of variables (\ref{cov}) imply
$$\int_{(-1,1)^c} H_x(dt)\; =\;  \car  \int_{-1}^x (x-z)^{\a-1} \hat{\varphi}(z)dz\; +\; \cah  \int_{x}^1 (z-x)^{\a-1} \hat{\varphi}(z)dz \; = \; 1.$$
On the other hand, in the case $\a >1,$ (\ref{cov}) and Lemma 2 show that
\begin{eqnarray*}
\int_{-1}^1 H_x^*(dt) &= &\cah \int_{-1}^1 (x-z)^{\a-1} \varphi(z)\,dz\; -\; \kappa_{\alpha,\rho}^*(x)\int_{-1}^1 \hat{\varphi}(y)\,dy\\
&= &1 \; -\;\lpa\frac{\Ga(1-\a\rho)2^{1-\alpha}}{\Ga(\a\hr)\Ga(1-\a)} \;+\;   (\a-1)\,\cah\int_{-1}^1 \hat{\varphi}(y)dy\rpa \int_1^x \psi_{\a,\rho}(t) dt\;=\;1,
\end{eqnarray*}
because
$$\int_{-1}^1 \hat{\varphi}(y)dy \;=\; 2^{1-\alpha}\, {\rm B}(1-\a\rho, 1-\a\hr) \;=\;  \frac{ 2^{1-\alpha} \Gamma(1-\a\rho)}{(1-\a)\cah \Gamma(\a\hr)\Gamma(1-\a)}\cdot$$
In the case $\a =1$, the measure $H_x^*$ has also total mass one by continuity. In the case $\a < 1,$ we find
$$\int_{-1}^1 H_x^*(dt)\; =\; 1 - \pb_x[T^* = \infty] \; =\; 1 \; -\;\frac{\Ga(1-\a\rho)2^{1-\alpha}}{\Ga(\a\hr)\Ga(1-\a)}\,\int_1^x \psi_{\a,\rho}(t) dt,$$
in accordance with Corollary 2 of \cite{BGR} and Corollary 1.2 of \cite{KPW}.

\end{remark}

\section{Proof of Theorem B}

\subsection{Proof of Part (a)} It is enough to consider the case $y >x,$ the case $x >y$ following from Hunt's switching identity - see e.g. Theorem II.5 in \cite{B}. By the Skorokhod continuity argument, it is also enough to consider $\a \neq 1.$ Reasoning as above, the D\'esir\'e Andr\'e equation yields
\begin{eqnarray*}
g(x,y) & = & c_\a\lpa \car\, (y-x)^{\a -1}\; - \; \int_{(-1,1)^c} u(t,y) H_x(dt)\rpa\\
&= &c_\a\lpa \car (y-x)^{\a -1} -  \car \int_{-\infty}^{-1} (y-t)^{\a -1} H_x(dt)  -  \cah\int_{1}^{+\infty} (t-y)^{\a -1} H_x(dt)\rpa
\end{eqnarray*}
with $c_\a = \Ga(1-\a).$ Changing the variable as in (\ref{cov}), we deduce 
\begin{eqnarray*}
g(x,y) & = & \Ga(1-\a)\, \car\,(y-x)^{\a -1} \lpa 1\,-\, \cah \int_{-1}^1 (z+t)^{\alpha-1}  \hat{\varphi}(t) dt\right)\\
&= & \Ga(1-\a)\,\car\,(y-x)^{\a -1} \lpa 1\,- \,\cah \int_{-1}^1 (z-s)^{\alpha-1} \varphi(s) ds\rpa
\end{eqnarray*}
and the result follows from Lemma 2, since $z > 1.$ 

\qed

\subsection{Proof of Part (b) in the case $\a < 1$} 

\subsubsection{The case $y>1$} Hunt's switching identity shows again that it is enough to consider the case $y > x.$ As above, the D\'esir\'e Andr\'e equation and the change of variable (\ref{cov}) imply 
\begin{eqnarray*}
g^*(x,y) & = & \Ga(1-\a) \,\car \lpa(y-x)^{\a -1}\, -\, \int_{-1}^{1} (y-t)^{\a -1} H_x^*(dt) \rpa\\
& = & \Ga(1-\a) \,\car \,(y-x)^{\a -1}\lpa 1 \, -\, \cah\int_{-1}^1 (z-u)^{\alpha-1} \varphi(u) du \rpa
\end{eqnarray*}
with $z > x > 1,$ and we can conclude by Lemma \ref{lem:phi}. 

\subsubsection{The case $y<-1$} Still using (\ref{cov}), we now have  
\begin{eqnarray*}
g^*(x,y) & = & \Ga(1-\a) \,\cah \lpa(x-y)^{\a -1}\,-\, \int_{-1}^{1} (t-y)^{\a -1} H_x^*(dt) \rpa\\
& = &  \Ga(1-\a) \,\cah \,(x-y)^{\a -1}\lpa 1\, -\, \int_{-1}^1 (z-u)^{\alpha-1} \varphi(u) du\rpa 
\end{eqnarray*} 
with $z \in (1,x),$ and we again conclude by Lemma \ref{lem:phi}. 
\qed

\subsection{Proof of Part (b) in the case $\a>1$.} We only consider the case $y > x.$ The case $x > y > 1$ is obtained by Hunt's switching identity and the case $y < -1$ by analogous computations. Proceeding as for Equation (\ref{Abel2}), we first deduce
$$g^*(x,y)\;=\;  \Ga(1-\a)\, \car\, \lpa(y-x)^{\a -1} - \int_{-1}^{1} (y-t)^{\a -1} H_x^*(dt) \rpa \, -\, \Gamma(1-\alpha)\kappa_{\alpha, \rho}^*(x).$$
Using Theorem 1 and the computations of the case $\alpha<1$, the expression may be transformed into
$$\frac{1}{\Ga(\a\rho)\Ga(\a\hr)} \left(\frac{y-x}{2}\right)^{\alpha-1} \int_1^{z} \psi_{\a,\rho} (t)\, dt\; -\; \Ga(1-\a)\,\kappa_{\alpha, \rho}^*(x) \lpa 1\, -\, \car \int_{-1}^1 (y-t)^{\a -1}\hat{\varphi}(t)dt\rpa.$$
The result follows from the hat version of Lemma 2  and the expression of $\kappa_{\alpha, \rho}^*(x).$

\qed

\section{Proof of the Corollaries}

\subsection{Proof of Corollary 1} By duality, it is enough to consider the case $x > y.$ From Part (a) of Theorem B and a change of variable, we see that $g(x,y)$ extends by continuity on the diagonal, with
$$g(y,y)\; =\; \frac{1}{(\a-1)\Ga(\a\rho)\Ga(\a\hr)} \lpa \frac{1-y^2}{2}\rpa^{\a -1}.$$
Moreover, it is clear that $g$ vanishes on the boundary $\{\vert x\vert =1\}\cup\{ \vert y\vert =1\}$ and is hence bounded on $(-1,1)\times(-1,1).$ By Proposition VI.4.11, Exercise VI.4.18 and Formula V.3.16 in \cite{BG}, we deduce
$$\pb_x[T_y < T]\; =\; \frac{g(x,y)}{g(y,y)}$$
and the conclusion follows by Theorem B.

\qed

\begin{remark} (a) In the case $\a \le 1,$ the process $L$ does not hit points, so that the problem is irrelevant. In general, one can ask for an evaluation of the probability
$\pb_x[ T_I < T]$ where $I$ is a closed subinterval of $(-1,1)$ not containing $x,$ and $T_I$ is its first hitting time. In the transient case $\a < 1,$ this problem is solved theoretically as a particular instance of the so-called condenser problem - see Formula (3.4) in \cite{CK}. It is an interesting open problem to find out an explicit formula in the real stable framework.\\

(b) By the Markov property, we can write down the following expression for the harmonic measure $H^{\{y\}}_x(dt)$ of the set $\{y\}\,\cup\,[-1,1]^c:$
\begin{equation}
\label{Harmy}
H^{\{y\}}_x(dt)\; =\; \rho(x,y)(\delta_{\{y\}}(dt) - H_y(dt))\; +\; H_x(dt).
\end{equation}
In particular, for every $x,y\in (-1,1),$ one has 
$$\pb_x [L_T \in dt, \, T < T_y]\; =\; H_x(dt) -\rho(x,y) H_y(dt).$$ 

\medskip

(c) It is interesting to mention that using the Gauss formula, we can deduce the asymptotic behaviour of $\pb_x[T_y > T]$ when $x\to y,$ which is fractional. For instance, if $y = 0,$ one has
$$\pb_x[T_0 > T]\, \mathop{\sim}_{x\rightarrow 0+}\, \frac{\Ga(2-\a)\Ga(\a\rho)}{\Ga(1-\a\hr)}\, (2x)^{\a-1}\quad \mbox{and}\quad \pb_x[T_0 > T]\, \mathop{\sim}_{x\rightarrow 0-}\, \frac{\Ga(2-\a)\Ga(\a\hr)}{\Ga(1-\a\rho)}\, \vert 2x\vert^{\a-1}.$$

\medskip

(d) By (\ref{Sca}) and spatial homogeneity, it is easy to deduce from Corollary 1 the following expression of ${\tilde \rho}(x,y) = \pb_x [T_y < \tau]$ where $\tau = \inf\{t > 0, \, L_t > 1\}:$ one finds
$${\tilde \rho}(x,y) \; =\; (\a-1) \lva \frac{x-y}{1-y}\rva^{\a -1} \int_0^{\lva\frac{1-x}{x-y}\rva} t^{\a\rho -1} (t+1)^{\a\hr -1}\, dt$$
if $x > y,$ and ${\tilde \rho}(x,y) = {\hat {\tilde \rho}}(-x,-y)$ if $x < y.$ When $y= 0,$ this is Theorem 1.5 in \cite{KPW}, correcting a misprint (the $1-1/x$ in the second integral should be $-1/x$) therein. Notice that Corollary 1.6 in \cite{KPW} is also analogously recovered from (\ref{Harmy}).

\end{remark}

\medskip

\subsection{Proof of Corollary 2} By the general theory of Martin boundary - see e.g. Theorem 1 in \cite{KuW}, we need to compute the Martin kernels
$$M_1(x)\; =\; \lim_{y \to 1} \,\frac{g(x,y)}{g(0,y)}\qquad \mbox{and}\qquad M_{-1} (x) \; =\; \lim_{y \to -1} \,\frac{g(x,y)}{g(0,y)}\cdot$$ 
Part (a) of Theorem B and a straightforward asymptotic analysis show that these Martin kernels exist and equal respectively
$$M_1(x)\; =\; (1-x)^{\a\rho -1}(1+x)^{\a\hr}\qquad\mbox{and}\qquad  M_{-1}(x)\; =\; (1+x)^{\a\hr -1}(1-x)^{\a\rho},$$
whence the result. \qed

\section{Final remarks}

In this section, we briefly describe the analogues of the above results in the case of semi-finite intervals and in the spectrally one-sided situation.

\subsection{The case of semi-finite intervals} By scaling and spatial homogeneity, one can deduce from Theorem A - either its Part (a) or its Part (b) - the following expression of the density of $L_{\tau}$ under $\pb_{x},$ where $x <1$ and $\tau =\inf\{ t >0, \, L_t > 1\}.$ One finds 
$$f_{L_{\tau}} (y)\, =\,\frac{\car(1-x)^{\a\rho}}{(y-1)^{\a\rho}(y-x)}\cdot$$
This expression has been found by several authors and can be obtained in different ways (see Exercise VIII.3 in \cite{B} and the references therein). Observe that it serves as a starting formula in \cite{Ro} in order to prove Part (a) of Theorem A. Notice last that the expression extends to the case with no negative jumps, by the Skorokhod continuity argument. In the relevant case with no positive jumps $\a > 1, \rho = 1/\a,$ the law of $L_\tau$ is a Dirac mass at one.

The Green function is
$$g_\tau(x,y) \; =\; \frac{(y-x)^{\a-1}}{\Ga(\a\rho)\Ga(\a\hr)}\, \int_0^{\frac{1-y}{y-x}} \psi_{\a,\rho}(t)\, dt$$
if $x < y < 1$ and $g_\tau(x,y) = {\hat g}_\tau(y,x)$ if $y < x <1.$ In the case $\a > 1,$ the analogue of Corollary 1 which is already given in Remark 6 (d) above, can then be recovered. Finally, one finds that the non-negative harmonic functions vanishing on $(1,+\infty)$ are of the type 
$$\lbd (1-x)^{\a\rho} \,+\, \mu(1-x)^{\a\rho -1}$$ 
with $\lbd, \mu \ge 0,$ in accordance with Theorem 4 in \cite{Si} and the paragraph thereafter.
   
\subsection{The case of stable processes with one-sided jumps}  By duality, it is enough to consider the two cases $\a < 1, \rho =1$ and $\a > 1, \rho = 1/\a.$

\subsubsection{The case $\a < 1, \rho =1$} It follows readily from the above paragraph that
$$h(x,y) \; =\; \frac{c_{\a,1}(1-x)^{\a}}{(y-1)^{\a}(y-x)}\,\Un_{\{y >1\}}$$
for all $x\in (-1,1).$  See also Example 3 in \cite{IW} and the references therein for the  expression of the density of $(L_{T-}, L_T)$ under $\pb_x.$ Similarly, one has  
$$h^*(x,y) \; =\; \frac{c_{\a,1}\vert 1+x\vert^{\a}}{(1+y)^{\a}(y-x)}\,\Un_{\{ |y|<1 \}}$$
for all $ x < -1$ and $h^*(x,y) = 0$ for all $x > 1.$ In accordance with the fact that $L$ is a subordinator, the Green function is
$$g(x,y) \; =\; \frac{(y-x)^{\a-1}}{\Ga(\a)} \, \Un_{\{x<y\}}$$
for all $x, y \in (-1,1),$ 
$$g^*(x,y) \; =\; \frac{(y-x)^{\a-1}}{\Ga(\a)} \lpa \Un_{\{x<y<-1\}}\, +\, c_{\a,1} \lpa\int_0^{\frac{\vert 1+x\vert(y-1)}{2}} \!\!\! t^{\a-1}(1+t)^{-1}\, dt \rpa\Un_{\{y>1\}}\rpa$$
for all $x <-1,$ and $g^*(x,y) = g(x,y)$ for all $x >1.$ The problem of Corollary 1 is irrelevant. Finally, the non-negative harmonic functions on $(-1,1)$ vanishing on $[-1,1]^c$ are constant multiples of $(1-x)^{\a-1}.$ 

\subsubsection{The case $\a > 1, \rho =1/\a$} Using Skorokhod continuity in Theorem A (a) and the absence of positive jumps, one has
$$H_x(dy)\; =\;c_{\a, 1-1/\a}(1-x)(1+x)^{\a-1}(1-y)^{-1}\vert y+1\vert^{1-\a}(x-y)^{-1}\Un_{\{y <-1\}}\, dy\; +\; \pb_x[ T_1 < T] \, \delta_1(dy).$$
Either taking the limit in Remark 6 (d) or integrating the first term, we can compute the weight of the Dirac mass, and find
$$H_x(dy)\; =\; c_{\a, 1-1/\a}(1-x)(1+x)^{\a-1}(1-y)^{-1}\vert y+1\vert^{1-\a}(x-y)^{-1}\Un_{\{y <-1\}}\, dy\; +\; \lpa\frac{x+1}{2}\rpa^{\a-1} \!\!\! \delta_1(dy).$$
The corresponding Green function is
$$g(x,y)\; =\; \frac{1}{\Ga(\a)}\lpa \lpa \frac{(1-y)(1+x)}{2}\rpa^{\a-1}\! -\; (x-y)^{\a-1}\Un_{\{x > y\}}\rpa.$$
The hitting probabilities are
$$\pb_x[T_y < T]\; =\; \lpa\frac{1+x}{1+y}\rpa^{\a-1}$$
for every $x \le y,$ which is also a consequence of a well-known result on scale functions - see e.g. Theorem VII.8 in \cite{B}, and$$\pb_x[T_y < T]\; =\; \lpa\frac{1+x}{1+y}\rpa^{\a-1} \!\!-\; \lpa\frac{2(x-y)}{1-y^2}\rpa^{\a-1}$$ 
for every $x > y.$ Finally, the non-negative harmonic functions on $(-1,1)$ which vanish on $[-1,1]^c$ are of the type $\lbd (1-x)^{\a-1}(1+x)^{\a -2} + \mu (1+x)^{\a-1}$
with $\lbd, \mu \ge 0.$\\

It is clear that $H_x^*(dy) =\delta_{-1}(dy)$ for all $x < -1.$ To compute $H_x^*(dy)$ for $x >1,$ let us introduce $\tau^* =\inf\{ t >0, \, L_t <1\}.$ The absence of positive jumps and the formula for semi-finite intervals imply after some computation
\begin{eqnarray*}
H_x^*(dy) &= & \Un_{\{ |y|<1\}} \Pb_x[L_{\tau^*} \in dy]  \;+\; \Pb_x[L_{\tau^*} <-1] \delta_{-1}(dy)\\
&= & c_{\a,1-1/\a}\!\lpa \frac{(x-1)^{\a-1} (1-y)^{1-\a} }{x-y}\,\Un_{\{|y|<1\}} dy + \lpa \int_{0}^{\frac{x-1}{x+1}} \!\!\! z^{\a-2} (1-z)^{1-\a}  dz \rpa\! \delta_{-1}(dy)\rpa,
\end{eqnarray*}
in accordance with Remark 3 in \cite{P1} - see also Proposition 1.3 in \cite{KPW}.

\bigskip

\noindent
{\bf Acknowledgements.}  Nous savons gr\'e \`a Jean Jacod d'un chat instructif sur la distance de Skorokhod. C. P. a b\'en\'efici\'e du support de la Chaire {\em March\'es en Mutation}, F\'ed\'eration Bancaire Fran\c{c}aise. Travail d\'edi\'e \`a l'association Laplace-Gauss.

\end{document}